\input amstex
\documentstyle{amsppt}
\magnification=\magstep1
 \hsize 13cm \vsize 18.35cm \pageno=1
\loadbold \loadmsam
    \loadmsbm
    \UseAMSsymbols
\topmatter
\NoRunningHeads
\title
Note on the Euler numbers and polynomials
\endtitle
\author
  Taekyun Kim
\endauthor

\abstract In this paper we investigate the properties of the Euler
functions. By using the Fourier transform for the Euler function,
we derive the interesting formula related to the infinite series.
Finally we give some interesting identities between the Euler
numbers and the second kind stirling numbers.

\endabstract
\thanks  2000 AMS Subject Classification: 11B68, 11S80
\newline
\endthanks
\endtopmatter

\document

{\bf\centerline {\S 1. Introduction/Definition}}
 \vskip 20pt
The constants $E_k$ in the Taylor series expansion
$$F(t)=\frac{2}{e^t+1}=\sum_{n=0}^{\infty}E_n\frac{t^n}{n!},
\text{ for $|t|\leq \pi$}, \text{ (see [1-31])},\tag1$$ are known as the Euler numbers.
The first few are $1,$ $-\frac{1}{2}$, $0$, $\frac{1}{4}, \cdots,$
and $E_{2k}=0$ for $k=1, 2, 3, \cdots.$ These numbers arise in the
series expansions of trigonometric functions, and are extremely
important in number theory and analysis. The Euler polynomials,
$E_n(x)$, are defined as
$$F(x,
t)=F(t)e^{xt}=\frac{2}{e^t+1}e^{xt}=\sum_{n=0}^{\infty}E_n(x)\frac{t^n}{n!},
\text{ for $x\in \Bbb R$} .\tag2$$  For $x\in\Bbb R$ with $0\leq
x<1$, the Euler polynomials are called the Euler functions. From
(1) and (2) we can derive
$$E_n(x)=\sum_{l=0}^n\binom{n}{l}E_lx^{n-l}, \text{ where $\binom {n}{l}=\frac{n(n-1)\cdots(n-l+1)}{l!}$.}
\tag3$$ Thus, we obtain  the distribution relation for the Euler
polynomials  as follows. For $d\in\Bbb N$ with $d\equiv 1 (\mod
2)$, we have
$$\sum_{k=0}^{d-1}(-1)^kE_{n}(\frac{x+k}{d})=d^{-n}E_n(x).$$
By (1), it is easy to see that the recurrence relation for the
Euler numbers is given by
$$E_0=1, \text{ and, }
\sum_{l=0}^n\binom{n}{l}E_l+E_n=2\delta_{0,n} \text{ where
$\delta_{0,n}$ is Kronecker symbol.}\tag4$$ From (3) and (4), we
note that
$$E_n(1)=\sum_{l=0}^n\binom{n}{l}E_l=-E_n, \text{ for $n\geq 1$}.
\tag5$$ Thus, we obtain the following lemma.

\proclaim{Lemma 1} For $n\in\Bbb N$, we have $E_n(1)=-E_n.$
\endproclaim
From (3) we can easily derive
$$\aligned
&\frac{dE_n(x)}{dx}=\frac{d}{dx}\sum_{k=0}^n\binom{n}{k}E_k
x^{n-k}=\sum_{k=0}^n\binom{n}{k}(n-k)E_kx^{n-k-1}\\
&=n\sum_{k=0}^n\frac{(n-1)!}{(n-k-1)!k!}E_kx^{n-1-k}=n\sum_{k=0}^{n-1}\binom{n-1}{k}E_kx^{n-1-k}
=nE_{n-1}(x).
\endaligned\tag6$$
By (6), we obtain the following proposition.

\proclaim{Proposition 2} For $n\geq 0$, we have
$$ \int_0^x E_{n}(t)dt=\frac{1}{n+1}E_{n+1}(x). \tag7$$
\endproclaim

In this paper we investigate the properties of the Euler
functions. By using the Fourier transform for the Euler function,
we derive the interesting formula related to the infinite series.
Finally we give some interesting identities between the Euler
numbers and the second kind stirling numbers.
 \vskip 20pt

{\bf\centerline {\S 2. Euler Functions }}\vskip 10pt

In  this section, we assume that $E_n(x)$ is the Euler function.
Let us consider the Fourier transform  for the Euler function,
$E_n(x),$ as follow. For $m\in \Bbb N$, the Fourier transform on
the Euler function is given by
$$E_m(x)=\sum_{n=-\infty}^{\infty}a_n^{(m)}e^{(2n+1)\pi ix}, \text{
($a_n^{(m)}\in\Bbb C$),}\tag8$$ where
$$a_n^{(m)}=\int_0^1 E_m(x)e^{-(2n+1)\pi ix}dx. \tag9$$
From (9), we note that
$$\aligned
&a_n^{(m)}=\int_0^1 E_m(x)e^{-\pi i (2n+1)x}dx\\
&=\left[\frac{E_{m+1}(x)}{m+1}e^{-\pi i(2n+1)x}\right]_0^1
+\frac{(2n+1)\pi i}{m+1}\int_0^1 E_{m+1}(x)e^{-(2n+1)\pi ix}dx\\
&=\frac{(2n+1)\pi i}{m+1}\int_0^1E_{m+1}(x)e^{-(2n+1)\pi ix}dx
=\frac{(2n+1)\pi i}{m+1}a_n^{(m+1)}.
\endaligned\tag10$$
Thus, we have
$$a_n^{(m)}=\frac{m}{(2n+1)\pi i}a_n^{(m-1)}=\frac{m(m-1)}{\left((2n+1)\pi i
\right)^2}a_n^{(m-2)}=\cdots=\frac{m!}{\left((2n+1)\pi i
\right)^{m-1}}a_n^{(1)}. \tag11$$

It is easy to see that
$$\aligned
&a_n^{(1)}=\int_0^1 E_1(x)e^{-(2n+1)\pi
ix}dx=\int_0^1(x-\frac{1}{2})e^{-(2n+1)\pi ix}dx\\
&=-\frac{1}{(2n+1)\pi i}\left[(x-\frac{1}{2})e^{-(2n+1)\pi
ix}\right]_0^1 +\frac{1}{(2n+1)\pi i}\int_0^1 e^{-(2n+1)\pi i x
}dx\\
&=\frac{2}{\left((2n+1)\pi i \right)^2} .
\endaligned\tag12$$
From (11) and (13), we can derive
$$a_n^{(m)}=2\frac{m!}{\left((2n+1)\pi i \right)^{m+1}}, \text{
$m\in\Bbb N$, and } a_n^{(0)}=\frac{2}{(2n+1)\pi i}. \tag13 $$

By (8) and (13), we see that
$$E_m(x)=m!2\sum_{n=-\infty}^{\infty}\frac{e^{(2n+1)\pi
ix}}{\left((2n+1)\pi i \right)^{m+1}}, \text{ for $ 0 \leq x
<1$}.$$

Therefore, we obtain the following theorem.

\proclaim{Theorem 3} For $ m \in \Bbb Z_{+}=\Bbb N \bigcup \{0\}$,
$x\in\Bbb R$ with $0\leq x <1$,  we have
$$E_m(x)=m!2\sum_{n=-\infty}^{\infty}\frac{e^{(2n+1)\pi
ix}}{\left((2n+1)\pi i \right)^{m+1}}.$$
\endproclaim

If we take $x=1$, then we have
$$E_m(1)=-m!2\sum_{n=-\infty}^{\infty}\frac{1}{\left((2n+1)\pi i \right)^{m+1}}.\tag14$$
By(14) and Lemma 1, we obtain the following  corollary.

\proclaim{Corollary 4} For $ m \in \Bbb Z_{+}=\Bbb N \bigcup
\{0\}$, we have
$$E_m=m!2\sum_{n=-\infty}^{\infty}\frac{1}{\left((2n+1)\pi i \right)^{m+1}}.$$
\endproclaim

From Corollary 4, we note that
$$E_{2m+1}=(-1)^{m+1}2 \frac{(2m+1)!}{\pi^{^{2m+2}}}\sum_{n=-\infty}^{\infty}\frac{1}
{(2n+1)^{2m+2}}.$$

Thus, we have
$$\sum_{n=-\infty}^{\infty}\frac{1}
{(2n+1)^{2m+2}}=(-1)^{m+1}
\frac{E_{2m+1}}{2(2m+1)!}\pi^{^{2m+2}}.\tag15$$

By (15), we obtain the following corollary.

\proclaim{Corollary 5} For $ m \in \Bbb Z_{+}$, we have
$$\sum_{n=1}^{\infty}\frac{1}
{(2n+1)^{2m+2}}=(-1)^{m+1}
\frac{E_{2m+1}}{4(2m+1)!}\pi^{^{2m+2}}.$$
\endproclaim

Note that
$$\aligned
&\frac{1}{1+e^{-x}}=\sum_{n=0}^{\infty}e^{-nx}(-1)^n=\sum_{n=0}^{\infty}(e^{-x})^n(-1)^n
=\sum_{n=0}^{\infty}\left(\sum_{k=0}^{\infty}\frac{(-1)^k}{k!}x^k\right)^n(-1)^n\\
&=\sum_{n=0}^{\infty}(-1)^n\left(\sum_{a_1+a_2+\cdots=n}\frac{n!}{a_1!
a_2!\cdots}\frac{(-1)^{a_1+2a_2+\cdots}
}{(1!)^{a_1}(2!)^{a_2}\cdots} \right)x^{a_1+2a_2+\cdots}.
\endaligned\tag16$$
Let $p(i,j): a_1+2a_2+\cdots=i, a_1+a_2+\cdots=j.$ From (16), we note that
$$\aligned
  &\frac{1}{1+e^{-x}}=\sum_{m=0}^{\infty}\left(\sum_{n=0}^m(-1)^n\sum_{p(m,n)}
  \frac{n!}{a_1!a_2!\cdots a_m!}\frac{(-1)^mx^m}{(1!)^{a_1} \cdots
  (m!)^{a_m}}\right)\\
  &=\sum_{m=0}^{\infty}(-1)^m\left(\sum_{n=0}^mn!(-1)^n\sum_{p(m,n)}\frac{m!}{a_1!a_2!\cdots
  a_m!}\frac{(-1)^m}{(1!)^{a_1}(2!)^{a_2}\cdots(m!)^{a_m}}\right)\frac{x^m}{m!}\\
  &= \sum_{m=0}^{\infty}(-1)^m\sum_{n=0}^m
  n!(-1)^ns_2(m,n)\frac{x^m}{m!},
\endaligned\tag17$$
where $s_2(m,n)$ is the second kind  stirling number.

By the definition of Euler number, we easily see that
$$\frac{1}{1+e^{-x}}=\frac{1}{2}\left(\frac{2}{1+e^{-x}}\right)=\frac{1}{2}\sum_{m=0}^{\infty}(-1)^mE_m\frac{x^m}{m!},
\text{ (see [5] )}.
\tag18$$

By comparing the coefficients of $\frac{t^n}{n!}$ on the both sides of (17) and (18), we obtain
$$E_m=2\sum_{n=0}^{\infty}(-1)^n n!s_2(m,n),$$
where $s_2(m,n)$ is the second kind stirling number.

\vskip 20pt

 \quad Taekyun Kim

   \quad
       Division of General Education-Mathematics,

\quad Kwangwoon University, Seoul, S. Korea \quad

\quad e-mail:\text{ tkkim$\@$kw.ac.kr; tkim64$\@$hanmail.net} \vskip10pt

 \enddocument